# Random environment on coloured trees

MIKHAIL MENSHIKOV,[1] DIMITRI PETRITIS[2] and STANISLAV VOLKOV[3]

[1]*Department of Mathematics, University of Durham, Durham DH1 3LE, U.K.*
*E-mail: Mikhail.Menshikov@durham.ac.uk*

[2]*Institut de Recherche Mathématique de Rennes, Université de Rennes 1, 35042, Rennes Cedex, France. E-mail: Dimitri.Petritis@univ-rennes1.fr*

[3]*Department of Mathematics, University of Bristol, BS8 1TW, U.K.*
*E-mail: S.Volkov@bristol.ac.uk*

In this paper, we study a regular rooted coloured tree with random labels assigned to its edges, where the distribution of the label assigned to an edge depends on the colours of its endpoints. We obtain some new results relevant to this model and also show how our model generalizes many other probabilistic models, including random walk in random environment on trees, recursive distributional equations and multi-type branching random walk on $\mathbb{R}$.

## 1. Introduction

Random walks in random environment have been studied for a long time (see Solomon [13] for such a random walk on $\mathbb{Z}$). One of the most natural extensions of this model is to consider a random walk in random environment on a tree (see, e.g., Lyons and Pemantle [9]). It turns out that the question of recurrence vs. transience of the walk is equivalent to the question of infiniteness vs. finiteness of certain sums of random variables assigned to the edges of the tree. In [9], it is assumed that all these random variables are i.i.d., which may be a fairly restrictive condition. Indeed, in the classical setup, the probability of jump from a given vertex $v$ through a certain edge is set to be equal to the ratio between the value assigned to this edge and the sum of the values assigned to all the edges adjacent to $v$. The assumption that these values are independent results in substantial restrictions on the possible random jump distributions assigned to vertices, in particular, the symmetry of such a distribution.

Our initial motivation in writing this paper was to overcome this difficulty. Additionally, we have managed to study the situations where the distribution of values assigned to an edge depends on its direction and even on the direction of the immediately preceding edge. This has resulted in the establishment of two phase transitions in the model, which, in turn, is useful for a variety of applications (not only random walks in random environment, as outlined in Section 5).

Formally, let $b \geq 2$ and consider a $b$-ary regular rooted tree $T = T_b$ with root $v_0$ and vertex set $\mathbb{V}$ (i.e., a tree in which all vertices have degree $b+1$, with the exception of the root, which has degree $b$). For any two vertices $v, u \in \mathbb{V}$, let $d(u,v)$ denote the distance between these two vertices, that is, the number of edges on the shortest path connecting $v$ and $u$. Let $\mathbb{V}_n$ denote the set of $b^n$ vertices at graph-theoretical distance $n$ from the root and write $|v| = d(v, v_0) = n$ when $v \in \mathbb{V}_n$. If two vertices $v$ and $w$ are connected by an edge, we write $v \sim w$ and let $\ell(v)$ denote the sequence of vertices of the unique self-avoiding path connecting $v$ to the root.

With the exception of the root, colour each vertex in one of $b$ distinct colours, from left to right, such that for every fixed vertex, each of its children has a different colour. For definiteness, colour the root in any of the $b$ colours. Denote the colour of vertex $v$ by $c(v) \in \{1, 2, \ldots, b\}$.

We are given $b^2$ positive-valued random variables of known joint distribution, which we denote $\bar{\xi}_{ij}$, $i, j = 1, 2, \ldots, b$. Now, to each unoriented edge $(u, v) \equiv (v, u)$ assign a random variable $\xi_{uv}$, such that

- for any edge $(u, v)$, where $u$ is *the parent* of vertex $v$, we have $\xi_{uv} \stackrel{\mathcal{D}}{=} \bar{\xi}_{c(u)\,c(v)}$;
- for any collection of edges $(u_1, v_1), (u_2, v_2), \ldots, (u_m, v_m)$ such that $u_i$ is the parent of $v_i$ for all $i$ and $u_i \neq u_j$ for all $i \neq j$, the random variables $\{\xi_{u_i v_i}\}_{i=1}^m$ are independent.

Here and throughout the paper, $X \stackrel{\mathcal{D}}{=} Y$ means that random variables $X$ and $Y$ have the same distribution. Note that we allow dependence between sibling edges.

For any $v \in \mathbb{V}$ and any $w \in \ell(v)$, let $\xi[w, v]$ equal the product of the random variables assigned to the edges of the subpath connecting $w$ to $v$. By convention, set $\xi[v, v] = 1$ and also denote $\xi[v] = \xi[v_0, v]$.

In this paper, we will answer the following two questions.

***Question 1.*** *When is* $\mathsf{Y} := \sum_{v \in \mathbb{V}} \xi[v]$ *finite a.s.?*

***Question 2.*** *Let $a > 0$. When is $\mathsf{Z} = \mathsf{Z}(a) := \mathrm{card}\{v \in \mathbb{V} : \xi[v] > a\}$ finite a.s.?*

The answers, of course, depend on the distribution of $\bar{\xi}_{ij}$'s and they are presented in Section 2, with the proofs given in Section 4, while Section 3 contains some auxiliary statements.

The study of the sum $\mathsf{Y}$ was a very important ingredient in the analysis carried out in [9] and it is essential for the investigation of random walks in random environment on trees, edge-reinforced random walks on trees and some other problems. At the same time, the quantity $\mathsf{Z}$ is relevant to first-passage percolation and branching random walks. Some other relevant models are also mentioned in [11].

In our paper, we also consider various applications and they are presented in Section 5.

## 2. Main results

We will first introduce an alternative colouring procedure, which is equivalent to the one described above in terms of the questions we ask, yet which has some advantages. Suppose that the colouring of the tree is done in a different manner than that described above. Namely, it is done recursively for $\mathbb{V}_1, \mathbb{V}_2, \ldots$, as follows. Suppose that the vertices up to level $n-1$ are already coloured. Next, colour the vertices of $\mathbb{V}_n$ randomly in such a way that whenever two vertices $v \in \mathbb{V}_n$ and $u \in \mathbb{V}_n$ share a common parent, they must have different colours. The distribution of the colouring is independent of the previous levels and is uniform, that is, we assign each of the allowed $(b!)^{b^{n-1}}$ colourings with equal probability. This process can be extended to infinity, thus producing the colouring of all vertices $v \in \mathbb{V}$.

Again, assign to each edge $(u,v)$ a random variable $\zeta_{uv}$ such that *conditioned on the colouring of the tree*, $\zeta_{uv}$ satisfies the two conditions on $\xi_{uv}$'s mentioned in the previous section. Similarly, compute $\zeta[w,v]$'s and $\zeta[v]$'s. Then, by construction, it is clear (e.g., by using coupling arguments) that for each $n$, the distribution of the *unordered set* $\{\zeta[v], \ v \in \mathbb{V}_n\}$ is the same as the distribution of $\{\xi[v], \ v \in \mathbb{V}_n\}$. Therefore, the answers to the questions above will be the same as in the original model. At the same time, the new model, which uses randomized colouring, has a significant advantage, namely that

$$\text{for any two } v, w \in \mathbb{V}_n, \qquad \zeta[v] \stackrel{\mathcal{D}}{=} \zeta[w], \tag{2.1}$$

though $\zeta[v]$ and $\zeta[w]$ are, of course, dependent. Thus, we will hereafter only work with the new, randomly coloured model. The probability $\mathbb{P}$ and the expectation $\mathbb{E}$ below will be with respect to the measure generated by a random colouring $\mathbf{c} = \{c(v), \ v \in \mathbb{V}\}$ and a random environment $\zeta = \{\zeta_{vw}, \ v,w \in \mathbb{V}, \ v \sim w\}$.

For $s \in [0, \infty)$, let

$$m(s) = \begin{pmatrix} \mathbb{E}(\bar{\xi}_{11})^s & \ldots & \mathbb{E}(\bar{\xi}_{1b})^s \\ \vdots & \ddots & \vdots \\ \mathbb{E}(\bar{\xi}_{b1})^s & \ldots & \mathbb{E}(\bar{\xi}_{bb})^s \end{pmatrix}$$

and let $\rho(s)$ be the largest eigenvalue of $m(s)$, which is positive by the Perron–Frobenius theorem, since all elements of the matrix are strictly positive and it is hence irreducible.

We will need the following regularity conditions. Let

$$\mathbb{D} = \{s \in \mathbb{R} : \mathbb{E}\bar{\xi}_{ij}^s < \infty \ \forall i, j \in \{1, 2, \ldots, b\}\}$$

and suppose that

$$\begin{aligned} &[0,1] \subseteq \mathbb{D}, \\ &0 \in \mathsf{Int}(\mathbb{D}), \\ &\mathbb{E}|\log \bar{\xi}_{ij}| < \infty \qquad \forall i,j \in \{1,2,\ldots,b\}, \\ &\mathbb{E}|\bar{\xi}_{ij} \log \bar{\xi}_{ij}| < \infty \qquad \forall i,j \in \{1,2,\ldots,b\}. \end{aligned} \tag{2.2}$$

We are now ready to present our main results.

**Theorem 1.** *Let* $Y = \sum_{v \in \mathbb{V}} \zeta[v]$ *and* $\lambda_1 = \inf_{s \in [0,1]} \rho(s)$.

(a) *If* $\lambda_1 < 1$, *then* $Y < \infty$ *a.s.*
(b) *If* $\lambda_1 > 1$ *and the conditions (2.2) in the next section are fulfilled, then* $Y = \infty$ *a.s.*

**Theorem 2.** *Let* $x > 0$, $Z(x) = \operatorname{card}\{v \in \mathbb{V} : \zeta[v] > x\}$ *and* $\lambda = \inf_{s \in [0,\infty)} \rho(s)$. *Additionally, suppose that conditions (2.2) are fulfilled.*

(a) *If* $\lambda < 1$, *then* $Z(x) < \infty$ *a.s.*
(b) *If* $\lambda > 1$, *then* $Z(x) = \infty$ *a.s.*

Note that we do not attempt to analyze here the situation in the critical case $\lambda = 1$ (resp. $\lambda_1 = 1$) The reason is that unlike the one-dimensional situation, the analysis becomes much more difficult here and we could not find any reasonable and interesting conditions which would ensure infiniteness of $Z(x)$ or $Y$.

## 3. Large deviations results

Let $v \in V_n$ and suppose that $\ell(v) = \{v_0, v_1, \ldots, v_{n-1}, v_n = v\}$. The random variables $c(v_i)$, $i = 1, 2, \ldots, n$, are then i.i.d. random variables with uniform distribution on the set $\{1, 2, \ldots, b\}$.

**Lemma 1.** *Let* $S_n = \sum_{i=1}^n \log(\zeta_{v_{i-1} v_i})$ *and*

$$k_n(s) = (\mathbb{E}(e^{sS_n}))^{1/n} = \left( \mathbb{E} \prod_{i=1}^n \zeta_{v_{i-1} v_i}^s \right)^{1/n}.$$

*Suppose that conditions (2.2) are fulfilled. Then,*

(a) $k(s) = \lim k_n(s) \in [0, \infty]$ *exists for all* $s$;
(b) $\Lambda(s) = \log \rho(s) - \log b = \log k(s) \in (-\infty, +\infty]$ *is convex;*
(c) *the rate function* $\Lambda^*(z) = \sup_{s \in \mathbb{D}} (sz - \Lambda(s))$, $z \in \mathbb{R}$, *is convex, lower semicontinuous and differentiable in* $\operatorname{Int}(\mathbb{D})$. *Moreover,*

$$\Lambda^*(z) = \begin{cases} s_0(z) z - \Lambda(s_0(z)), & \text{if } z \geq \Lambda'(0), \\ 0, & \text{if } z \leq \Lambda'(0), \end{cases}$$

*where* $s_0(z)$ *is the solution of the equation* $z - \Lambda'(s) = 0$;
(d) *for all* $a > 0$,

$$\lim_{n \to \infty} \frac{1}{n} \log \mathbb{P}\left( \frac{S_n}{n} \geq \log a \right) = -\Lambda^*(\log a).$$

**Remark 1.** The statement of Lemma 1 holds simultaneously for all possible colourings of the root.

**Proof of Lemma 1.** First, for any $s_1, s_2 \in \mathbb{D}$ with $s_1 < s_2$, the segment $[s_1, s_2]$ belongs to $\mathbb{D}$, hence for any $\alpha \in (0,1)$, we have $k_n^n(\alpha s_1 + (1-\alpha)s_2) = \mathbb{E}(\exp(\alpha s_1 S_N))\exp((1-\alpha)s_2 S_N)) \leq [\mathbb{E}(\exp(s_1 S_N))]^\alpha [\mathbb{E}(\exp(s_2 S_N))]^{1-\alpha}$ from which logarithmic convexity of $k_n$ follows.

Suppose $c(v_0) = \alpha$. Then, $k_n(s) = \frac{1}{b}(e_\alpha^T m(s)^n e)^{1/n}$, where

$$e_\alpha = \begin{pmatrix} 0 \\ \vdots \\ 0 \\ 1 \\ 0 \\ \vdots \\ 0 \end{pmatrix} \leftarrow \alpha\text{th position}, \qquad e = \begin{pmatrix} 1 \\ 1 \\ \vdots \\ 1 \end{pmatrix}. \tag{3.3}$$

Now, $m(s) = \rho_1(s)P_1 + \rho_2(s)P_2 + \ldots$, where $(\rho_i)$ are the eigenvalues of $m(s)$, ordered so that $\rho_1 > |\rho_2| \geq |\rho_3| \geq \ldots$, and $P_i$ denotes the projection on the $i$th eigenspace corresponding to the $i$th eigenvalue $\rho_i$. Notice that $\rho \equiv \rho_1 > 0$, the image space of $P_1$ is 1-dimensional and since $m(s)_{ij} > 0$ for all $i,j$, we have

$$(P_1 e)_i > 0 \quad \text{for all } i \text{ and } |\rho_2| < \rho_1.$$

Hence, $k(s) = \frac{1}{b}\rho(s) \in (0, \infty]$ for all $s$ and is log-convex, as the limit of a log-convex function is log-convex.

Finally, (d) follows from the Gärtner–Ellis theorem (see, e.g., Lemma V.4, page 53 in den Hollander [6]) under conditions (2.2). $\square$

Note that we can rewrite $\Lambda^*$ as

$$\Lambda^*(z) = \sup_{s \geq 0}[sz - \log(\rho(s)/b)]$$

$$= -\log \inf_{s \geq 0} \frac{\rho(s)\mathrm{e}^{-sz}}{b}. \tag{3.4}$$

Recall that

$$\lambda = \inf_{s \geq 0} \rho(s),$$

$$\lambda_1 = \inf_{s \in [0,1]} \rho(s) \geq \lambda.$$

**Corollary 1.** *For any $\tilde{\lambda} < \lambda_1$ and $\alpha \in \{1,\ldots,b\}$, there exists a $y \in (0,1]$ and a positive integer $n$ such that for any $v, u \in \mathbb{V}$ such that $u \in \ell(v)$, $c(u) = \alpha$ and $d(u,v) = n$, where*

$$\mathbb{P}(\zeta[u,v] > y^n) \geq \frac{\tilde{\lambda}}{(by)^n}.$$

**Proof.** Lemma 1 yields that for any small $\varepsilon > 0$ and all $y > 0$, there exists an $n_0 = n_0(\varepsilon, y)$ such that, for every $n \geq n_0$,

$$e^{n\varepsilon} \mathbb{P}(S_n/n \geq \log y) \geq \exp\{-n\Lambda^*(\log y)\}$$
$$= \exp\left\{n \inf_{s \geq 0}[\log(\rho(s)/b) - s\log y]\right\}$$
$$= \left(\inf_{s \geq 0} \frac{\rho(s)y^{-s}}{b}\right)^n = \frac{1}{(yb)^n}\left(\inf_{s \geq 0} \rho(s)y^{1-s}\right)^n,$$

whence, for any $v, u \in \mathbb{V}$ such that $u \in \ell(v)$ and $|v| = |u| + n$,

$$\mathbb{P}(\zeta[u,v] \geq y^n) \geq \left[\frac{e^{-\varepsilon}}{yb} \times \inf_{s \geq 0} \rho(s)y^{1-s}\right]^n.$$

Now, since $\log \rho(s)$ is convex, it follows from the proof of the lemma on page 129 in Lyons and Pemantle [9] that

$$\max_{0 < y \leq 1} \inf_{s \geq 0} \rho(s)y^{1-s} = \min_{0 \leq s \leq 1} \rho(s) = \lambda_1.$$

Consequently, by choosing $y \in (0,1]$ at the point where this maximum is achieved, and $\varepsilon > 0$ very small, we ensure that for all large $n$,

$$\mathbb{P}(\zeta[u,v] \geq y^n) \geq \frac{\tilde{\lambda}}{(yb)^n}. \tag{3.5}$$

□

## 4. Proofs of the main theorems

**Proof of Theorem 1.** (a) Suppose that $\lambda_1 < 1$. We can then fix an $s \in (0,1)$ such that $\rho(s) < 1$. Suppose that the root has colour $\alpha$. Then,

$$\sum_{v \in \mathbb{V}_n} \zeta^s[v] = \sum_{\ell(v) = (v_0,\ldots,v_n):\ v = v_n \in \mathbb{V}_n} \zeta^s_{v_0 v_1} \zeta^s_{v_1 v_2} \cdots \zeta^s_{v_{n-1} v_n}$$

and hence, by construction of the colouring of the tree, we have

$$\mathbb{E}\left(\sum_{v\in\mathbb{V}_n}\zeta^s[v]\right) = \sum_{\ell(v)=(v_0,\ldots,v_n):\ v=v_n\in\mathbb{V}_n} m_{\alpha c(v_1)}(s)m_{c(v_1)c(v_2)}(s)\cdots m_{c(v_{n-1})c(v_n)}(s),$$

$$= e_\alpha^T m^n(s)e,$$

where $e_\alpha$ and $e$ are defined in (3.3). Now, since $\rho(s) < 1$, $\sum_{n=1}^\infty m^n(s) < \infty$, therefore, by Fubini's theorem, $\mathbb{E}(\sum_v \zeta^s[v]) < \infty$ and hence, $\sum_v \zeta^s[v]$ is finite a.s. This implies that $\zeta[v] \geq 1$ for *only finitely many* $v$'s and therefore, there exists an $N$ such that whenever $v \notin \mathbb{V}_1 \cup \cdots \cup \mathbb{V}_N$, it follows that $\zeta[v] < 1$, hence $\zeta^s[v] > \zeta[v]$. Consequently,

$$\mathsf{Y} = \sum_{i=1}^N \left(\sum_{v\in\mathbb{V}_i}\zeta[v]\right) + \sum_{i=N+1}^\infty \left(\sum_{v\in\mathbb{V}_i}\zeta[v]\right)$$

$$< \sum_{i=1}^N \left(\sum_{v\in\mathbb{V}_i}\zeta[v]\right) + \sum_{i=N+1}^\infty \left(\sum_{v\in\mathbb{V}_i}\zeta^s[v]\right) < \infty,$$

where the last inequality follows from the fact that $\sum_v \zeta^s[v]$ is finite a.s.

(b) Since $\lambda_1 > 1$, by Corollary 1, there exist an $\varepsilon > 0$, $0 < y \leq 1$ and an $n$ such that for any $v, u \in \mathbb{V}$ satisfying $u \in \ell(v)$ and $|v| = |u| + n$,

$$\mathbb{P}(L[u,v]) \geq \frac{1+\varepsilon}{(by)^n} =: q,$$

where

$$L[u,v] := \{\zeta[u,v] \geq y^n\}.$$

Let us construct an embedded branching process, with members of generation $j$ denoted $M_j$, as follows. The root of the tree $v_0$ is the sole member of generation 0, that is, $M_0 = \{v_0\}$. For $j \geq 1$, let

$$M_j = \{u \in \mathbb{V}_{jn} : \exists w \in M_{j-1} \text{ such that } \mathbb{V}_{(j-1)n} \cap \ell(u) = \{w\} \quad (4.6)$$
$$\text{and } L[w,u] \text{ occurs}\}.$$

The process $|M_j|$ can be minorized by an *independent* branching process with uniformly bounded number of descendants whose average is equal to

$$\mu := b^n \times q = (1+\varepsilon)y^{-n} > 1,$$

which is a supercritical process surviving with a positive probability, say, $p_S > 0$. Moreover,

$$\left\{\lim_{j\to\infty}\frac{|M_j|}{\mu^j} > 0\right\} = \{\text{the process survives}\} \quad \text{a.s.}$$

by the Kesten–Stigum theorem (see, e.g., [1], page 192). This, in turn, implies that there is a positive $\delta > 0$ such that, with probability at least $p_S/2 > 0$ for all $j$ sufficiently large, we have $|M_j| \geq \delta \mu^j$. Consequently, since for each $v \in M_j$, we have $\zeta[v] \geq (y^n)^j$, it follows that

$$\mathsf{Y} \geq \sum_{v \in M_j} \zeta[v] \geq \delta \mu^j (y^n)^j = \delta(1+\varepsilon)^j \to \infty \qquad \text{as } j \to \infty,$$

with positive probability. Now, the set $\{\mathsf{Y} = \infty\}$ is a tail event and the random variables at different generations are independent, hence its probability satisfies the 0–1 law and we obtain the required result. $\square$

**Proof of Theorem 2.** (a) Recall that for any $v \in \mathbb{V}_n$, the quantity $p_n = \mathbb{P}(\zeta[v] > x)$ does not depend on $v$ and observe that

$$\mathbb{E}\mathsf{Z}(x) = \sum_{n=1}^{\infty} b^n p_n.$$

Since $\lambda < 1$, from (3.4), we have that for a small $z < 0$ and a very small $\varepsilon > 0$,

$$-\Lambda^*(z) < \log \frac{1 - 2\varepsilon}{b}.$$

Set $y = e^z < 1$ and apply Lemma 1(d) to obtain that for all large $n$,

$$\frac{1}{n} \log \mathbb{P}(\zeta[v] \geq y^n) \leq \log \frac{1 - \varepsilon}{b}.$$

This yields

$$b^n \mathbb{P}(\zeta[v] \geq y^n) \leq (1 - \varepsilon)^n$$

and since $x > 0$ and $y < 1$ implies that $p_n = \mathbb{P}(\xi[v] > x) \leq \mathbb{P}(\xi[v] \geq y^n)$ for large $n$, we have $\mathbb{E}\mathsf{Z}(x) < \infty$, so $\mathsf{Z}(x) < \infty$ a.s.

(b) Now, since $\lambda > 1$, from (3.4), we have that for a small $z > 0$ and a very small $\varepsilon > 0$,

$$-\Lambda^*(z) > \log \frac{1 + 2\varepsilon}{b}.$$

As before, we set $y = e^z > 1$ and apply Lemma 1(d) to obtain that there exists an $n$, which we shall now fix, such that

$$\frac{1}{n} \log \mathbb{P}(\zeta[v] \geq y^n) \geq \log \frac{1 + \varepsilon}{b} \quad \Longrightarrow \quad b^n \mathbb{P}(\zeta[v] \geq y^n) \geq (1 + \varepsilon)^n. \qquad (4.7)$$

Next, we construct a branching process that is almost identical to the one constructed in the proof of Theorem 1. Again, provided that $u \in \ell(v)$ and $|v| - |u| = n$, we introduce the event

$$L[u, v] := \{\zeta[u, v] \geq y^n\}, \qquad (4.8)$$

whose probability is at least $(1+\varepsilon)^n/b^n$, according to (4.7). Let the root of the tree $v_0$ be the unique member of generation 0, that is, $M_0 = \{v_0\}$. Similarly to the previous proof, for $j \geq 1$, let $M_j$ be defined by (4.6). Then, the process $|M_j|$ can again be minorized by a supercritical independent branching process, with average number of descendants equal to $\mu := (1+\varepsilon)^n > 1$, which survives with positive probability $p_S > 0$. On the event $\sum_{j=1}^{\infty} |M_j| = \infty$ of survival, for any $x > 0$, there exists $j_0 = j_0(x)$ such that for all $j \geq j_0$, we have $v \in M_j$, implying that $\zeta[v] \geq y^{nj} > x$. Consequently,

$$\mathbb{P}(\mathsf{Z}(x) = \infty \text{ for all } x > 0) \geq p_S > 0.$$

Taking into account the fact that the event $\{\mathsf{Z}(x) = \infty \text{ for all } x > 0\}$ is a tail event and variables at different generations are independent, we conclude that for any $x > 0$,

$$\mathbb{P}(\mathsf{Z}(x) = \infty) = 1.$$

□

## 5. Applications

Here, we show how Theorems 1 and 2 can be applied to obtain some of the already known facts, as well as to establish new results in various applications of probability theory. Throughout this section, we will assume that the regularity conditions (2.2) are satisfied.

### 5.1. Random walk in random environment

Let $u$ be a vertex of the coloured tree $T$. For every $v \sim u$, define $p_{uv} \in (0,1)$ such that $\sum_{v:\, v \sim u} p_{uv} = 1$. For definiteness, denote the parent of $u$ by $u^*$ and the children of $u$ by $u^1, u^2, \ldots, u^b$ (also, when $u$ is the root $v_0$ of the tree, set $u^* \equiv u$). Now, suppose that for each $u$,

$$\mathbf{p}(u) = (p_{uu^*}, p_{uu^1}, p_{uu^2}, \ldots, p_{uu^b}) \in (0,1)^{b+1}$$

is a $(b+1)$-dimensional random variable. Obviously, the set of components of $\mathbf{p}(u)$ is *dependent*, since they must sum to 1.

Suppose that the distribution of $\mathbf{p}(u)$ depends only on the colour $c(u)$ of the vertex $u$. Additionally, suppose that the random variables $\{\mathbf{p}(u), u \in \mathbb{V}\}$ are independent.

Now, define a random walk $X(k)$ in a random environment on the coloured tree $T$ by letting $X(0) = v_0$ and

$$\mathbb{P}(X(k+1) = v \mid X(k) = u) = \begin{cases} p_{uv}, & \text{if } u \sim v, \\ 0, & \text{otherwise,} \end{cases}$$

where we set $v_0 \sim v_0^* = v_0$. This model is similar to that of Lyons and Pemantle [9]; however, we do not require as much independence or symmetry for the distribution of the jumps to children as is required in [9]. Additionally, we allow jump distributions to depend on the type of the vertex. On the other hand, in [9], more general trees are considered, while we restrict ourselves to regular trees.

We want to establish when the walk in the random environment is transient (resp. recurrent). For $i = 1, 2, \ldots, b$, let

$$(\bar{\xi}_{i1}, \bar{\xi}_{i2}, \ldots, \bar{\xi}_{ib}) \stackrel{\mathcal{D}}{=} \left( \frac{p_{uu^1}}{p_{uu^*}}, \frac{p_{uu^2}}{p_{uu^*}}, \ldots, \frac{p_{uu^b}}{p_{uu^*}} \right)$$

whenever $c(u) = i$. Also, let $m(s)$, $\rho(s)$ and $\lambda_1$ be the same as defined in Section 2.

**Proposition 1.** *The random walk in random environment described above is a.s. positive recurrent when $\lambda_1 < 1$ and a.s. transient when $\lambda_1 > 1$.*

**Proof.** We will use the standard electric network representation of the random walk by replacing each edge of the coloured tree $T$ with a resistor, such that their conductances satisfy the formula

$$\frac{C_{uu^*}}{C_{uu^i}} = \frac{p_{uu^*}}{p_{uu^i}},$$

where, again, $u^1, \ldots, u^b$ are the children of vertex $u$ and $u^*$ is its parent (see [7]). These equations are satisfied when, for any $u \in \mathbb{V}_n$, $n \geq 1$, with $\ell(u) = \{u_0 \equiv v_0, u_1, u_2, \ldots, u_{n-1}, u_n \equiv u\}$, we have

$$C_{u_{n-1}u_n} = \prod_{i=0}^{n-1} \frac{p_{u_i u_{i+1}}}{p_{u_i u_{i-1}}}.$$

Note that when $u_i = v_0$, we need to set $p_{u_i u_{i-1}} = p_{v_0 v_0}$.

Now, to each edge $(u^*, u)$ where $u^*$ is the parent of $u$, assign a random variable with distribution $\bar{\xi}_{c(u^*)c(u)}$. Then, $C_{u_{n-1}u_n}$ is equal to the product of the random variables assigned to edges of the path connecting $v_0$ to $u_n = u$. Theorem 1 implies that whenever $\lambda_1 < 1$, $\mathsf{Y} = C := \sum_{x,y} C_{x,y} < \infty$ a.s. and then there exists a stationary probability measure $\pi$ such that $\pi_x = C_x/C$, where

$$C_x = \sum_{y: \, y \sim x} C_{xy}.$$

Therefore, the random walk is positive recurrent.

The reverse statement (concerning transience for $\lambda_1 > 1$) follows from a slight modification of the proof of part (i) of Theorem 1 of [9], effectively using the estimate (3.5), since transience is equivalent to establishing finiteness of the effective resistance $R_{\text{eff}}$. □

*Example.* Consider a random walk in a random environment on a coloured binary tree ($b = 2$). Suppose that from a vertex of type 1, the walk always goes down with probability $\frac{1}{2}$ and up with probabilities $\frac{1}{4}$ to either of its children. From a vertex of type 2, the walk goes up and right with probability $\frac{1}{4}$, down with probability $\frac{3}{4}\eta_v$ and up and left with probability $\frac{3}{4}(1-\eta_v)$, where $\eta_v$ are i.i.d. random variables distributed uniformly on $[h,1]$, $h \in (0,1)$. Then,

$$m(s) = \left( \mathbb{E}\left(\frac{1-\eta}{\eta}\right)^s \quad \mathbb{E}\left(\frac{1}{3\eta}\right)^s \right)^{2^{-s}}.$$

It is easy to verify that if $\rho(s)$ is the largest eigenvalue of $m(s)$, then $\lambda_1 = \inf_{s \in [0,1]} \rho(s)$ is smaller than 1 whenever $h > h_{cr} = 0.417\ldots$, thus the walk is positive recurrent for almost every environment when $h > h_{cr}$ and transient for almost every environment when $h < h_{cr}$.

## 5.2. Recursive distributional equations

It turns out that our construction on randomly coloured trees may be used to answer the question about the existence of solutions of certain functional equations, which will be described below.

Let $\Xi$ be the $b \times b$ matrix of random variables $\bar{\xi}_{ij}$, $i,j = 1,2,\ldots,b$. We want to find a $b$-dimensional random vector $Y = (Y_1,\ldots,Y_b)^T$, independent of $\Xi$, such that

$$1 + \sum_{j=1}^{b} \bar{\xi}_{ij} Y_j \stackrel{\mathcal{D}}{=} Y_i \qquad \text{for} \quad i = 1, 2, \ldots, b, \tag{5.9}$$

which can be expressed in vector form as

$$e + \Xi Y \stackrel{\mathcal{D}}{=} Y,$$

where $e$ has been defined by (3.3). Equation (5.9) is a special case of a more general recursive distributional equation which have been widely studied; for example, see [10], where one can find sufficient conditions for the existence of its solution (Theorem 4.1). At the same time, for equation (5.9), we essentially obtain a *criterion* for this existence.

Note that it will be essential that the tree is coloured; otherwise, we would have been able to solve (5.9) only in the one-dimensional case $b = 1$.

**Proposition 2.** *Let* $\mathsf{Y}$ *be the quantity defined in Question 1 of Section 1. Then, equation (5.9) has a solution if and only if* $\mathsf{Y} < \infty$ *a.s.*

**Proof.** First, suppose that $\mathsf{Y} < \infty$ a.s. For each $i = 1,\ldots,b$, let $Y_i = \sum_{v \in \mathbb{V}} \zeta[v]$ when $c(v_0) = i$ and suppose that different $Y_i$'s are constructed using independent random variables. By assumption, each $Y_i < \infty$ a.s. It is now easy to see that $1 + \sum_{i=1}^{b} \bar{\xi}_{ij} Y_j$ indeed has the distribution of $Y_i$, hence $Y = (Y_1,\ldots,Y_b)^T$ is a solution of (5.9).

Second, suppose that $\mathsf{Y} = \infty$ a.s. and suppose that there exists a solution $\hat{Y} = (Y_1, Y_2, \ldots, Y_n)^T$ of equation (5.9). Construct a $b$-ary tree with $c(v_0) = 1$ and assigned random variables $\bar\xi$'s, as described in the Introduction. Also, for each $n$ and each $v \in \mathbb{V}_n$, let $Q(v)$ be an independent random variable with the distribution of $Y_{c(v)}$ and denote

$$\mathsf{Y}_1^{(<n)} = \sum_{v \in \mathbb{V}_0 \cup \cdots \cup \mathbb{V}_{n-1}} \zeta[v],$$

$$\tilde{\mathsf{Y}}_1^{(<n)} = \sum_{v \in \mathbb{V}_0 \cup \cdots \cup \mathbb{V}_{n-1}} \zeta[v] + \sum_{v \in V_n} \zeta[v] Q(v),$$

for $n = 1, 2, \ldots$. Since $\hat{Y}$ is a solution to the problem, it follows by induction on $n$ that $\tilde{\mathsf{Y}}_1^{(<n)}$ must have the distribution of $Y_1$ and this is true for *all* $n$. Now, observe that $\mathsf{Y}_1^{(<n)} \leq \tilde{\mathsf{Y}}_1^{(<n)}$. At the same time, $\lim_{n \to \infty} \mathsf{Y}_1^{(<n)} = \mathsf{Y} = \infty$ a.s., by assumption. Hence, $\tilde{\mathsf{Y}}_1^{(<n)}$, which equals to $Y_1$ in distribution, is larger than a random variable equal to $\infty$ a.s., which is impossible. □

Let $\lambda_1$ be as defined in Section 2. Then, Theorem 1 yields the following.

**Corollary 2.** *If $\lambda_1 < 1$, then there exists a solution to equation (5.9). In contrast, if $\lambda_1 > 1$, then equation (5.9) has no solution.*

### 5.3. First-passage percolation

Here, we show how our techniques can extend the results of the first-passage percolation theory to the situation where one allows *negative passage times*. For each edge $(u, v)$ of the coloured tree $T$, where $u$ is the parent of $v$, let $\tau_{uv}$ denote the passage time from vertex $u$ to vertex $v$. Allow these times to be also negative, for example, indicating a "speeding up" of a walker. Suppose, for simplicity, that the $\tau_{uv}$'s are all independent, while their distribution depends on the colour of the endpoints, thus being one of the $b^2$ possible types. Let

$$R(t) = \left\{ u \in \mathbb{V} \colon \sum_{(v,w) \in \ell(u)} \tau_{vw} \leq t \right\}$$

be the set of vertices of the tree which are reachable in time $t$. The primary question is whether $R(t)$ is finite, since it can easily be infinite, due to the negative passage times.

To answer this, for all $u$ and $v$ such that $u \sim v$ and $u$ is the parent of $v$, set $\bar\xi_{c(u)c(v)} \stackrel{\mathcal{D}}{=} e^{-\tau_{uv}}$. The following statement is straightforward.

**Proposition 3.** *$R(t)$ is finite a.s. if and only if the quantity $\mathsf{Z}(e^{-t})$ defined in Question 2 of Section 1 is finite a.s.*

### 5.4. Multi-type branching random walks on $\mathbb{R}^1$

The literature on branching random walks is fairly extensive and a similar model to the one which follows was considered, for example, in [4], although somewhat different questions were investigated in that paper. A very similar model was also considered in Biggins [2, 3].

Suppose that we are given $b^2$ positive-valued random variables $\eta_{ij}$, $i,j = 1, 2, \ldots, b$, with known joint distribution. Consider a process on $\mathbb{R}$ which starts with a single particle of type $i \in \{1, 2, \ldots, b\}$ located at point $X^{(0)} = 0 \in \mathbb{R}$. The particle splits into $b$ new particles, one of each type $1, 2, \ldots, b$. If the positions of the new particles of types $1, 2, \ldots, b$ are denoted $X_1^{(1)}, X_2^{(1)}, \ldots, X_b^{(1)}$, respectively, then the distributions of jumps $X_j^{(1)} - X^{(0)}$ are independent for different $j$'s and have the distribution of $\eta_{ij}$. After this, each of the new particles behaves in the same way as the original particle, so by time $t \in \{1, 2, \ldots\}$, we will have exactly $b^t$ particles $X_1^{(t)}, \ldots, X_{b^t}^{(t)}$ located somewhere on $\mathbb{R}$.

Set $\bar{\xi}_{ij} = \exp(-\eta_{ij})$, $i, j = 1, \ldots, b$. The following statement is then obvious.

**Proposition 4.** *Suppose that* $\mathsf{Z}(1)$, *as defined in Question 2 of Section 1, is finite a.s. All particles will then eventually be on the positive semi-axis a.s., that is,*

$$\mathbb{P}\left(\exists N: \ \forall t \geq N \ \min_{i \in \{1,2,\ldots,b^t\}} X_i^{(t)} \geq 0\right) = 1.$$

In fact, we can strengthen this result. Let

$$\mu_t = \min_{i=1,\ldots,b^t} X_i^{(t)}$$

be the minimum displacement of our multi-type branching random walk (see, e.g., [5] and references therein). As before, set $\bar{\xi}_{ij} = \exp(-\eta_{ij})$, $i, j = 1, \ldots, b$, and let $\rho(s)$ be the largest eigenvalue of

$$m(s) = \begin{pmatrix} \mathbb{E}\mathrm{e}^{-s\eta_{11}} & \cdots & \mathbb{E}\mathrm{e}^{-s\eta_{1b}} \\ \vdots & \ddots & \vdots \\ \mathbb{E}\mathrm{e}^{-s\eta_{b1}} & \cdots & \mathbb{E}\mathrm{e}^{-s\eta_{bb}} \end{pmatrix}.$$

For $x \in \mathbb{R}$, let

$$\lambda^{(x)} = \inf_{s \geq 0} \mathrm{e}^{sx} \rho(s) \tag{5.10}$$

and observe that $\lambda^{(x)}$ is non-decreasing in $x$. Note that if the joint distribution of $\eta_{ij}$'s is not degenerate, then $\rho(s)$ is strictly log-convex in $s$ and therefore, there exists a unique $x_0$ such that $\lambda^{(x_0)} = 1$.

**Proposition 5.** *Under the non-degeneracy condition above,*

$$\lim_{t \to \infty} \frac{\mu_t}{t} = x_0 \qquad a.s.,$$

*where $x_0$ is the unique solution of the equation $\lambda^{(x_0)} = 1$.*

**Proof.** For each $x \in \mathbb{R}$, we can define a new multi-type branching random walk with step sizes equal to $\eta_{ij}^{(x)} = \eta_{ij} - x$ for every $i$ and $j$. We can then naturally couple the new walk $X_k^{(t;\,x)}$ with the original one by setting $X_k^{(t;\,x)} = X_k^{(t)} - tx$. Observe that the largest eigenvalue $\rho^{(x)}(s)$ of the matrix $m^{(x)}(s)$ for this modified walk, whose entries are $(\mathbb{E}e^{-s(\eta_{ij}-x)})_{i,j=1}^b$, equals $e^{sx}\rho(s)$, hence the value of $\lambda$ needed for Theorem 2 is given by (5.10).

Suppose $x < x_0$. Then, $\lambda^{(x)} < 1$, whence, by Theorem 2, $Z(1) < \infty$ a.s. and, by Proposition 4,

$$\mathbb{P}(\mu_t - tx \geq 0 \text{ for all sufficiently large } t) = 1. \tag{5.11}$$

To prove the complementary statement, we need to improve slightly the proof of part (b) of Theorem 2. First, choose $x > x_0$ yielding $\lambda^{(x)} > 1$ and replace the event (4.8) by

$$\tilde{L}[u,v] := \{\zeta[u,v] \geq y^n, \text{ and } \zeta[u,w] > \nu \text{ for all } w \in \ell(v) \text{ with } |w| > |u|\}.$$

One can choose the constants $\nu > 0$ and $\varepsilon > 0$ sufficiently small that $\mathbb{P}(\tilde{L}[u,v]) > (1+\varepsilon)^n/b^n$ still holds. We can then construct sets $M_j$ defined by (4.6) with $L[u,v]$ replaced by $\tilde{L}[u,v]$.

Let $j_0$ be sufficiently large that $\nu y^{nj_0} > 1$ (recall that $y > 1$). On the event of survival of the process $|M_j|$, for each $j$, we have $M_j \neq \emptyset$ and $M_{j+1} \neq \emptyset$, hence there exist a $u \in M_j \subseteq \mathbb{V}_{nj}$ and a $v \in M_{j+1} \subseteq \mathbb{V}_{n(j+1)}$ such that $u \in \ell(v)$ and $\tilde{L}[u,v]$ occurs. Consequently, for every $t \geq nj_0$ such that $nj \leq t < n(j+1)$, there is a $w \in \mathbb{V}_t$ such that $w \in \ell(v)$, whence $\zeta[w] = \zeta[u,w]\zeta[u] \geq \nu y^{nj} > 1$. On the other hand, $\zeta[w] > 1$ implies that $\mu_t - tx < 0$. Therefore, we have proved that the event

$$\{\mu_t - tx < 0 \text{ for all sufficiently large } t\} \tag{5.12}$$

has positive probability, since the branching process minorizing $|M_j|$ is supercritical. However, the event (5.12) is a tail event, so it must have probability 1. Together with (5.11), this completes the proof of Proposition 5. □

### 5.5. Number theory: $5x+1$ Collatz-type problem

Fix an odd positive integer $q$ and define the following map:

$$T_q : x \to \begin{cases} x/2, & \text{if } x \text{ is even}, \\ qx+1, & \text{if } x \text{ is odd}. \end{cases}$$

The famous, yet unresolved, Collatz problem (see, e.g., [8] for hundreds of references to papers and short descriptions of their content) states that if one sequentially applies mapping $T_3$ to any positive integer, then it will eventually arrive at the cycle $1 \to 4 \to$ $2 \to 1$. On the other hand, a similar mapping $T_5$ is conjectured to "explode", that is, for most positive integers $\lim_{n\to\infty} T_5^{(n)}(x) = \infty$ [12].

Another conjecture made in [14] states that the density of those numbers $x \in \mathbb{Z}_+$ for which $\lim_{n\to\infty} T_5^{(n)}(x) < \infty$ has a "Hausdorff dimension" of approximately 0.68. This conjecture was made based on a construction of a probabilistic "equivalent" of mapping $T_5$, leading to a special case of the model studied in connection to Question 2. For more details, see [14].

## Acknowledgements

The authors wish to thank Maury Bramson, Robin Pemantle and the anonymous referee for useful comments. D. Petritis acknowledges partial support from the European Science Foundation (grant 420 of the program "Phase Transitions and Fluctuation Phenomena for Random Dynamics in Spatially Extended Systems").